\documentclass[12pt,a4paper,English]{article}
\usepackage{latexsym}
\usepackage{hyperref}

\usepackage{amsfonts,amssymb}
\usepackage{amsmath}
\def\boldsymbol#1{\mbox{\boldmath{$#1$}}}

\newcommand{\qed}{\hspace*{5 mm}$\triangle$}

\newcommand{\FF}{\mathbb{F}}

\newtheorem{theorem}{Theorem}[section]

\newtheorem{lemma}[theorem]{Lemma}
\newcommand{\pr}{\indent{\bf Proof. \ }}
\newcommand{\la}{\langle}
\newcommand{\ra}{\rangle}
\newcommand{\bx}{{\bf x}}
\newcommand{\bv}{{\boldsymbol v}}

\newcommand{\bc}{{\boldsymbol c}}
\newcommand{\by}{{\boldsymbol y}}
\newcommand{\bp}{{\boldsymbol p}}
\newcommand{\bz}{{\boldsymbol z}}

\newcommand{\bbz}{{\boldsymbol 0}}
\newcommand{\bbo}{{\boldsymbol 1}}

\newcommand{\be}{{\boldsymbol e}}
\newcommand{\supp}{\mbox{\rm supp}}

\newcommand{\Stab}{\mbox{\rm Stab}}

\newcommand{\Orb}{\mbox{\rm Orb}}
\newcommand{\rank}{\mbox{\rm rank}}
\newcommand{\wt}{\mbox{\rm wt}}
\newcommand{\mA}{\mbox{\rm A}}
\newcommand{\mC}{\mbox{\rm C}}
\newcommand{\mS}{\mathcal{S}}

\title{On the number of Steiner triple systems $S(2^m-1,3,2)$ of rank $2^m - m + 2$ over $\FF_2$
\footnote{The research was carried out at the IITP RAS at the
expense of the Russian Foundation for Sciences (project No. 14-50-00150).}}

\author{D.Zinoviev}

\begin{document}
\maketitle

\begin{abstract}
We obtain the number of different Steiner triple systems $S(2^m-1,3,2)$ of rank
$2^m-m+2$ over the field $\FF_2$.
\end{abstract}

\section{Introduction}

A Steiner System $S(v,k,t)$ is a pair $(X,B)$ where $X$ is a
set of $v$ elements and $B$ is a collection of $k$-subsets
(blocks) of $X$ such that every $t$-subset of $X$ is contained
in exactly one block of $B$. A System $S(v,3,2)$ is called a Steiner
triple system (briefly STS$(v)$). Tonchev in his work \cite{Ton}
enumerated all different Steiner triple systems STS$(v)$ of order
$v= 2^m-1$ whose $2$-rank (i.e. rank over the field $\FF_2$),
is equal to $2^m - m$. The case of rank $2^m-m+1$ has been considered
in \cite{ZZ1}. In that paper all different Steiner triple systems STS$(v)$
of this rank were enumerated and moreover the exact number of different
triple systems was obtained \cite{ZZ2}. The next value of rank $2^m-m+2$
was treated in \cite{ZZ3}, where all Steiner triple systems of this rank
were enumerated and the number of different systems was estimated
from above. For more on Steiner systems see the survey on the
Steiner systems \cite{CD}. The goal of this work is to obtain the number of different
Steiner triple systems STS$(v)$ of order $v=2^m-1$ and rank $2^m-m+2$
over the field $\FF_2$.

\section{Preliminary results and terminology}

Let $E_q$ be an alphabet of size $q$: $$E_q=\{0,1, \ldots, q-1\},$$
in particular, $E=\{0,1\}$. A code $C$ of length $n$ over $E_q$ is
an arbitrary subset of $E_q^n$. Denote the $q$-ary code
$C$ of length $n$ with the minimal (Hamming) distance $d$ and
cardinality $N$ as an $(n, d, N)_q$-code. In particular, when $q=2$
denote it as an $(n, d, N)$-code. Denote by $\wt(\bx)$ the Hamming
weight of vector $\bx$ over $E_q$, and by $d(\bx, \by)$ the Hamming
distance between the vectors $\bx, \by \in E_q^n$.
For a binary code $C$ denote by ${\la}C{\ra}$ the linear envelope of words
of $C$ over the Galois Field $\FF_2$. The dimension of the space ${\la}C{\ra}$ is
called the {\em rank} of code the $C$ over $\FF_2$ (alternatively $2$-rank) and
is denoted by $\rank\,(C)$. Throughout this paper we will use "rank"
instead of "$2$-rank".

Denote by $(n,w,d,N)$ a binary constant weight code $C$ of length
$n$, with $N$ codewords of a fixed weight $w$ and the minimal distance $d$.

Let $J=\{1,2,\ldots,n\}$ be a set of coordinate positions $E_q^n$.
Denote by $\supp(\bv) \subseteq J$ the support of a vector
$\bv = (v_1,\ldots,v_n) \in E^n$, i.e. the set of indices of
its non-zero coordinate positions:
$$
\supp(\bv) = \{i:~v_i \neq 0\}.
$$

A binary $(n,d,N)$-code $C$, which is a linear $k$-dimensional space
over $\FF_2$, is denoted as $[n,k,d]$-code. Let
$$
({\bx \cdot \by}) = x_1y_1 + \cdots + x_ny_n
$$
be the scalar product over $\FF_2$ of the binary vectors
$\bx = (x_1,\ldots,x_n)$ and $\by = (y_1,\ldots,y_n)$.
For any (linear, non-linear or constant weight) code $C$ of length $n$ let
$C^{\perp}$ be its dual code:
$$
C^{\perp}~=~\{\bv \in \FF_2^n:~(\bv \cdot \bc) = 0, ~\forall ~\bc \in C\}.
$$
It is clear that $C^{\perp}$ is a $[n,n-k,d^{\perp}]$-code with the
minimal distance $d^{\perp}$, and where $k=\rank\,(C)$.

Let $\mS_n$ be the permutation group of $n$ elements $J_n=\{0,1,\ldots,n-1\}$ (thus $|\mS_n|=n!$).
For any $i\in J_n$ and $\pi\in\mS_n$ let $\pi(i)$ be the image of $i$ under the action of $\pi$.
Define the action of $\mS_n$ on $E^n$ as the permutation of coordinates, i.e.
for any $\bx = (x_0,\ldots,x_{n-1})\in E^n$, and any $\pi\in\mS_n$ let
$$
\pi(\bx) = (x_{\pi^{-1}(0)},x_{\pi^{-1}(1)},\ldots,x_{\pi^{-1}(n-1)}).
$$
For any set $X\in E^n$ and any $\pi\in\mS_n$, let
$$
\Stab(X) = \{ \pi\in\mS_n ~|~ \pi(X) = X\}, ~\mbox{where}~ \pi(X) = \{ \pi(x) ~|~ \mbox{for all} ~x\in X\}.
$$

\begin{lemma}\label{lem:dim2}
Let $E^3$ be a vector space of dimension three. Then it has seven subspaces of dimension two.
\end{lemma}

\pr
Suppose that $E^3 = {\la}\bx_1,\bx_2,\bx_3{\ra}$. Then the subspaces are given by
$$
{\la}\bx_1,\bx_2{\ra}, {\la}\bx_1,\bx_3{\ra}, {\la}\bx_2,\bx_3{\ra},{\la}\bx_1,\bx_2+\bx_3{\ra},
{\la}\bx_1+\bx_2,\bx_3{\ra},{\la}\bx_1+\bx_3,\bx_2{\ra},{\la}\bx_1+\bx_2,\bx_2+\bx_3{\ra}.
$$
\qed

In the present work, we deal with several classes of different
codes. To avoid extra notations, we use the same notations
$C$, $L$, $V$, $W$ or $S$, for the codes (i.e. for the sets of vectors)
as well as for the matrices, formed by the codewords of these codes,
i.e. where the codewords have some ordering.

Let $L$ be an MDR $(3,2,8^2)_8$-code over the alphabet $E_8$.
Note that the code $L$ is given by a Latin square of order 8.
Denote by $\Gamma_L$ the number of different squares.

Define the mapping $\varphi$ of elements $0,1,\ldots,7$ of $E_8$ into $E^8$:
$$
\varphi(0) = (1,0,\ldots,0),~~\varphi(1) = (0,1,0,\ldots,0),~~\ldots ,~~\varphi(7) = (0,\ldots,0,1)
$$
and extend it coordinatewise to the vectors $\bc = (c_1, \ldots, c_n)$ over $E_8$:~
$$
\varphi(\bc) = (\varphi(c_1), \ldots, \varphi(c_n)).
$$
For a given code $(3,2,64)_8$-code $L$, define the constant weight $(24,3,4,64)$-code,
denoted by $C(L)$, as follows:
$$
C(L)~=~\{\varphi(\bc):~\bc \in L\}.
$$
Every codeword $\bc$ of code $C(L)$, is split into blocks
of length $8$ so that $\bc = (\bc_1, \bc_2, \bc_3)$ where $\wt(\bc_i) = 1$ for $i = 1,2,3$.

Let $V=W_2^8$ be the set of words of length $8$ and weight $2$ (it is obvious that $|W_2^8|=28$).
Split it into seven (trivial) subcodes $P_i$,~$i=1, \ldots, 7$ (i.e. into seven parallel
classes of length $8$ and weight $2$) where each subcode has parameters $(8,2,4,4)$, i.e.
$$
V~=~\bigcup_{i=1}^7 \;P_i = \{P_1,P_2,\ldots,P_7\}.
$$
Write $V=(P_1,\ldots,P_7)$ when the order of the parallel classes is fixed.
For any permutation $\tau\in \mS_7$ of indices $\{1,2,\ldots,7\}$, let
\begin{equation}\label{eq:perm}
\tau\ast V ~=~ (P_{\tau^{-1}(1)},P_{\tau^{-1}(2)},\ldots,P_{\tau^{-1}(7)}).
\end{equation}

Define vectors
$$\be_1=(1,0,0,0,0,0,0,0),$$
$$\be_2=(0,1,0,0,0,0,0,0),$$
$$ \ldots $$
$$\be_7=(0,0,0,0,0,0,0,1).$$

Let $\Gamma_V$ be the number of different such partitions $V=\{P_1,\ldots,P_7\}$.
It is known (see, for example Lemma 5 of \cite{ZZ3}) that there exist $\Gamma_V = 6240$
different partitions which form $6$ orbits under the action of the permutation group $\mS_8$
of coordinate positions.
Let $V^{(1)},V^{(2)},\ldots,V^{(6)}$ be the non-equivalent partitions so that the size of their
$\mS_8$-orbits are respectively equal
$$
\{630,420,2520,30,1680,960\}.
$$
These partitions are (the first row corresponds to $V^{(1)}$, the second to $V^{(2)}$
and so on)
$$
\begin{array}{ccccccc}
 (81 42 24 18, & 82 41 28 14, & 84 48 22 11, & 88 44 21 12, & 90 60 0\mA 05, & \mA 0 50 09 06, & \mC 0 30 0\mC 03),\\
 (81 44 22 18, & 82 48 21 14, & 84 41 28 12, & 88 42 24 11, & 90 60 0\mA 05, & \mA 0 50 09 06, & \mC 0 30 0\mC 03),\\
 (81 44 28 12, & 82 48 21 14, & 84 41 22 18, & 88 42 24 11, & 90 60 09 06,   & \mA 0 50 0\mA 05, &\mC 0 30 0\mC 03),\\
 (81 42 28 14, & 82 41 24 18, & 84 48 22 11, & 88 44 21 12, & 90 60 0\mA 05, & \mA 0 50 09 06, &\mC 0 30 0\mC 03),\\
 (81 44 22 18, & 82 60 14 09, & 84 50 21 0\mA, & 88 42 24 11, & 90 41 28 06, & \mA 0 48 12 05, &\mC 0 30 0\mC 03),\\
 (81 60 18 06, & 82 44 30 09, & 84 41 28 12, & 88 50 22 05, & 90 42 21 0\mC, & \mA 0 48 14 03, &\mC 0 24 11 0\mA),
\end{array}
$$
where every parallel class is presented as an 8-digit number, two digits per word of length $8$.

\section{The Main result}

Suppose $S_v = S(v,3,2)$ is a Steiner triple system of order $v=2^m-1$ and
not greater than $2^m - m + 2$. This means that the dual code $S_v^{\perp}$ contains
a subcode $[v,m-3,d^{\perp}]$, which is denoted by $\mathcal{A}_m$ and whose minimal
distance $d^{\perp}$ is given by: $d^{\perp} = (v+1)/2 = 2^{m-1}$ [4].
More precisely, $\mathcal{A}_m$ contains the non-zero words of the same
weight $2^{m-1}$, i.e. the code is a subcode of a well known
linear equidistant Hadamard code. Thus, without loss of generality
we can assume that the code $\mathcal{A}_m$ is generated
by the following matrix:
\begin{equation}\label{eq:3.1}
\left[
\begin{array}{ccccccccc}
\bf{1}   &\bf{1}  &\bf{1}   &\bf{1}  &\ldots &\bf{0} &\bf{0} &\bf{0}  &\bf{0}\\
\bf{1}   &\bf{1}   &\bf{0}   &\bf{0}   &\ldots &\bf{1} &\bf{1} &\bf{0}  &\bf{0}\\
\ldots &\ldots &\ldots &\ldots &\ldots &\ldots &\ldots &\ldots &\ldots\\
\bf{1}   &\bf{0}  &\bf{1}   &\bf{0}  &\ldots &\bf{1} &\bf{0} &\bf{1}  &\bf{0}\\
\end{array}
\right],
\end{equation}
where $\bbo=(1,1,1,1,1,1,1,1)$ and $\bbz=(0,0,0,0,0,0,0,0)$.

\medskip
Theorem 3 of \cite{ZZ3} provides the estimate of the number $M_v$ of different
Steiner triple systems $S(v,3,2)$ of order $v=2^m-1=8u+7\ge 15$ of rank not greater than
$v-m+3$, whose orthogonal code $\mathcal{A}_m$ is of the form (\ref{eq:3.1}):
\begin{eqnarray*}
M_v&=&30\cdot\left(\Gamma_L\right)^k\times
\left(7!\cdot\Gamma_V\right)^u\\
&=&30\cdot\left(108776032459082956800\right)^{k} \times
\left(31449600\right)^{u},
\end{eqnarray*}
where $k=u(u-1)/6$.

In the present work, we obtain the exact formula for the number of different triple
systems of rank $v-m+3$. The similar formula for the rank $v-m+2$ was obtained in \cite{ZZ2},
and for the rank $v-m+1$ in \cite{Ton}. The main result of this paper is the following theorem:

\begin{theorem}\label{theo:1.0}
The number $M_{v,3}^{(o)}$ of different Steiner triple systems $S(v,3,2)$
of order $v=2^m-1=8u+7\geq 15$ and of rank $v-m+3$ is equal to
$$
M^{(o)}_{v}~=~ \frac{v! M_{v,3}} {(u(u-1)(u-2) \cdots (u+1)/2)\cdot (8!)^u \cdot 7!},
$$
where $M_{v,3}$ is the number of different Steiner triple systems
$S(v,3,2)$ of order $v$, of rank $v-m+3$ and orthogonal to the code
$\mathcal{A}_m$ of the form (\ref{eq:3.1}), and we have that
\begin{eqnarray*}
M_{v,3}
&=& 30 \times (7!\cdot 6240)^{u} \times (108776032459082956800)^{k} \\
&-& 210\times (4!\cdot 3!)^u \times 5040^u \times 576^{4k}\times 2^{u-m+3}\\
&+& 210\times 48^u\times 2100^u\times 16^{4k}\times 2^{2(u-m+3)}.
\end{eqnarray*}
where 
$k=u(u-1)/6$.
\end{theorem}

This theorem directly follows from Theorem \ref{theo:1.1}, whose proof
requires a few auxiliary statements.\\

Define the mapping of linear spaces
\begin{equation}\label{eq:psiE8}
\psi:~  E^8 \mapsto  E^8/{\la}\bbz,\bbo{\ra}.
\end{equation}
Thus any vector $\bx\in E^8$ is identified with $\bar{\bx}$,
where $\bar{\bx}=\bx +\bf{1}$, the coordinatewise addition in $E^8$.
For any vector $\bx = (\bx_1 ~|~ \ldots ~|~\bx_u ~|~ \bx_{u+1})\in E^v$
extend the mapping $\psi$ coordinatewise to the space $E^v$:
\begin{equation}\label{eq:psiEv}
\psi(\bx) =  (\psi(\bx_1) ~|~ \ldots ~|~\psi(\bx_u) ~|~ \bx_{u+1}).
\end{equation}

For an arbitrary vector $\bx\in E^{24}$ let $\mathcal{L}(\bx)$
be a subset of the set $\mathcal{L}$ of all different codes $L$ such that
$C(L)$ is orthogonal to the vector $\bx$. Let $\bx,\by,\bz\in E^{24}$ then define
$\mathcal{L}(\bx,\by) = \mathcal{L}(\bx)\cap \mathcal{L}(\by)$ and
$\mathcal{L}(\bx,\by,\bz) = \mathcal{L}(\bx)\cap \mathcal{L}(\by)\cap\mathcal{L}(\bz)$.

\begin{lemma}\label{lem:1.1}
Let $\bx=(\bx_1 |~\bx_2 |~\bx_3)$, $\by=(\by_1 |~\by_2 |~\by_3)$ and
$\bz=(\bz_1 |~\bz_2 |~\bz_3)$ be pairwise orthogonal arbitrary vectors
of length $24$ and weight $12$ split into blocks of length $8$ so
that the weight of every block $\bx_i,\by_i,\bz_i$ ($i=1,2,3$)
is equal to $4$. Suppose we have that $\psi(\bx_i),\psi(\by_i),\psi(\bz_i)$
are linearly independent in the space $E^8/{\la}\bbz,\bbo{\ra}$.
Suppose that
$$
\wt(\bx_i+\by_i) ~=~ \wt(\bx_i+\bz_i) ~=~ \wt(\by_i+\bz_i) ~=~ 2,
$$
for $i=1,2,3$. Then
\begin{itemize}
\item $|\mathcal{L}(\bx)| ~=~ |\mathcal{L}(\by)| ~=~ |\mathcal{L}(\bz)|=576^4$;
\item $|\mathcal{L}(\bx,\by)| ~=~ 16^4$, and moreover $\mathcal{L}(\bx,\by) ~=~
\mathcal{L}(\bx,\bx+\by) ~=~ \mathcal{L}(\bx+\by,\by)$;
\item $|\mathcal{L}(\bx,\by,\bz)|= 1$.
\end{itemize}
\end{lemma}

\pr
Since the orthogonality (of two vectors) is preserved by the permutation
of coordinates, without loss of generality we can assume that
\begin{equation}\label{eq:2.0}
\bx = (1,1,1,1,0,0,0,0~|~1,1,1,1,0,0,0,0~|~1,1,1,1,0,0,0,0).
\end{equation}

Suppose that the code $L$ (i.e. the Latin square) is presented
as a matrix $\{a_{i,j}\}$, $i,j=0,\ldots,7$.
Note that for the codeword $(i,j,a_{i,j})\in L$, the vector
$\phi((i,j,a_{i,j}))=(\phi(i),\phi(j),\phi(a_{i,j}))$
is orthogonal to the vector (\ref{eq:2.0})
if and only if either two symbols from the triple
$\{i,j,a_{i,j}\}$ belong to the set $\{0,1,2,3\}$,
or none of these symbols belong to $\{0,1,2,3\}$.

Consider the $i$-th row of the Latin square $L$ of order 8:
\begin{equation}\label{eq:2.1}
(a_{i,0},a_{i,1},a_{i,2},a_{i,3},a_{i,4},a_{0,5},a_{0,6},a_{0,7}).
\end{equation}
Then for $i=0,1,2,3$, we have that
$$
\{a_{i,0},a_{i,1},a_{i,2},a_{i,3}\}=\{4,5,6,7\} ~~\mbox{and}~~ \{a_{i,4},a_{i,5},a_{i,6},a_{i,7}\}=\{0,1,2,3\}.
$$
Similarly, for $i=4,5,6,7$, we have that
$$
\{a_{i,0},a_{i,1},a_{i,2},a_{i,3}\}=\{0,1,2,3\} ~~\mbox{and}~~ \{a_{i,4},a_{i,5},a_{i,6},a_{i,7}\}=\{4,5,6,7\}.
$$
Thus one can see that the square matrix $L$ of order 8 is split into 4 submatrices of order 4,
so that each of these submatrices is a Latin square of order 4. We obtain that $|\mathcal{L}(\bx)|=576^4$,
since $576$ is the number of different Latin squares of order $4$.
The remaining part of the proof coincides with that of Lemma 1 of \cite{ZZ2}.

The equality $\mathcal{L}(\bx,\by) ~=~ \mathcal{L}(\bx,\bx+\by) ~=~ \mathcal{L}(\bx+\by,\by)$
is obvious.
\qed\\

For an arbitrary vector $\bx\in E^{8}$ of weight $4$ and an arbitrary partition
$V=\{P_1,\ldots,P_7\}$ of $W_2^8$ we say that the parallel class is $\bx$-even (respectively
$\bx$-odd) if all words of the parallel class are orthogonal to the vector $\bx$
(respectively not orthogonal). Note that among the parallel classes $P_i$ there
exist at most four $\bx$-odd and three $\bx$-even classes.
For an arbitrary partition $V$ denote by
$$
\mathcal{X}(V) = \{\bx\in E^{8} ~|~ \wt(\bx) = 4\}
$$
the set of all vectors such that among the parallel classes $P_i$, $i=1,2,\ldots,7$ there exist
four $\bx$-odd and three $\bx$-even parallel classes.
The vector $\bp\in E^7$ of parities of the parallel classes is called the parity vector. Thus
$$
\bp = \bp\,(\bx,V) = (p_1,p_2,\ldots,p_7),
$$
where $p_i$ is $0$ if the corresponding parallel class $P_i$ is $\bx$-even
and is $1$ when $P_i$ is $\bx$-odd.
Obviously
\begin{equation}\label{eq:par1}
\mbox{if}~ \bx\in\mathcal{X}(V)~\mbox{then}~ \bar{\bx}\in\mathcal{X}(V)~ \mbox{and}~ \bp(\bx,V) = \bp(\bar{\bx},V).
\end{equation}

\begin{lemma}\label{lem:1.2}
For all 6 nonequivalent (under the action of the permutation group $\mS_8$) partitions $V^{(i)}$,
where $i=1,2,3,4,5,6$, of the set $W_2^8$, the cardinality of the set $|\mathcal{X}(V^{(i)})|$
takes the following values
$$
\{6,2,2,14,0,0\}.
$$
\end{lemma}

\pr
Direct computations.
\qed\\

One can check that adding the zero vector and the vector of ones
$\bbz,\bbo\in E^8$ turns the sets $\mathcal{X}(V^{(i)})$ ($i=1,2,3,4$)
into the vector spaces whose dimension is equal to $\{3,2,2,4\}$). Namely
$$
{\la}\mathcal{X}(V^{(i)}){\ra} = \mathcal{X}(V^{(i)}) ~\cup~ \{\bbz,\bbo\}.
$$
We also have that the dimensions of the spaces ${\la}\psi(\mathcal{X}(V^{(i)})){\ra}$ ($i=1,2,3,4$)
ere respectively equal to $\{2,1,1,3\}$.

Recall that the set $S^{(3)}$ is the Steiner triple system $S(7,3,2)$ (see [4, $\S 4$]).
There is one such system up to equivalence and $30$ different systems.
The default system is equal to
\begin{equation}\label{eq:S3}
\begin{array}{ccccccc}
(1 & 0 & 1 & 0 & 1 & 0 & 0)\\
(0 & 1 & 0 & 1 & 1 & 0 & 0)\\
(0 & 1 & 1 & 0 & 0 & 1 & 0)\\
(1 & 0 & 0 & 1 & 0 & 1 & 0)\\
(1 & 1 & 0 & 0 & 0 & 0 & 1)\\
(0 & 0 & 1 & 1 & 0 & 0 & 1)\\
(0 & 0 & 0 & 0 & 1 & 1 & 1).
\end{array}
\end{equation}

This system is boolean and the set $S^{(3)\perp} = (S^{(3)})^\perp$ (dual to $S^{(3)}$)
of non-zero vectors coincides with $S^{(3)}$ where every word is inverted (i.e. a vector of ones is added).
The following facts are well known and can be verified directly.

\begin{lemma}\label{lem:1.2A}
The linear space $S^{(3)\perp}\subset E^7$ has the following properties:
\begin{itemize}
\item $\Stab(S^{(3)}) = \Stab(S^{(3)\perp})\subset S_7$;
\item $|\,\Stab(S^{(3)\perp})| = 168$;
\item For any non-zero $\bx_1,\bx_2\in S^{(3)\perp}$ we have
$$
|\,\Stab(\la\bx_1,\bx_2\ra)| = 48;
$$
where $\la\bx_1,\bx_2\ra$ is a $2$-dimensional subspace of $S^{(3)\perp}$;
\item The action of the group $\Stab(S^{(3)\perp})$ on the set of $2$-dimensional
subspaces of $S^{(3)\perp}$ is transitive.
\end{itemize}
\end{lemma}

We also need some simple facts about partitions $V^{(1)}$ and $V^{(4)}$.
All of them can be directly verified.

\begin{lemma}\label{lem:1.3}
The partitions $V^{(1)}$ and $V^{(4)}$ have the following properties
\begin{itemize}
\item the parity vectors of $V^{(1)}$ are equal to
$$
\begin{array}{ccccccc}
(0 & 0 & 1 & 1 & 1 & 1 & 0)\\
(1 & 1 & 0 & 0 & 1 & 1 & 0).
\end{array}
$$
\item the parity vectors of $V^{(4)}$ are equal to
\[
\begin{array}{ccccccc}
(0 & 1 & 0 & 1 & 0 & 1 & 1)\\
(1 & 0 & 1 & 0 & 0 & 1 & 1)\\
(1 & 0 & 0 & 1 & 1 & 0 & 1)\\
(0 & 1 & 1 & 0 & 1 & 0 & 1)\\
(0 & 0 & 1 & 1 & 1 & 1 & 0)\\
(1 & 1 & 0 & 0 & 1 & 1 & 0)\\
(1 & 1 & 1 & 1 & 0 & 0 & 0).
\end{array}
\]
Note that this is equal to the non-zero vectors of $S^{(3)\perp}$.
\end{itemize}
\end{lemma}

\medskip

\begin{lemma}\label{lem:1.4}
For any partition $V$ and vector $\bx\in\mathcal{X}(V)$
\begin{itemize}
\item for any permutation $\pi\in\mS_8$ of coordinates, we have
$$
\bp\,(\pi(\bx),\pi(V)) ~=~ \bp\,(\bx,V);
$$
\item for any permutation $\tau\in\mS_7$, we have
$$
\bp\,(\bx,\tau\ast V) ~=~ \tau^{-1}(\bp\,(\bx,V));
$$
\end{itemize}
\end{lemma}

\pr
The first statement follows from the fact that for any permutation $\pi$
$$
(\pi(\bx)\cdot\pi(\by)) ~=~ (\bx\cdot\by).
$$
Thus, the parity of $\bx$ and a parallel class $P$ is the same as that of $\pi(\bx)$
and $\pi(P)$.

The second one follows from the definition of $\tau\ast V$ and the parity vector.
\qed
\medskip

For any set $X\subset E^n$ and vector $\by\in E^m$ set
$$
X \times \by ~=~ \{(\bx\,|\,\by) ~|~ \bx\in X\}.
$$

\begin{lemma}\label{lem:1.5}
Let $V = (P_1,P_2,\ldots,P_7)$ be a partition. Consider the set
$$
\mathcal{P} ~=~ \bigcup_{i=1}^7 P_{i}\times \be_i.
$$
For any $\bx\in\mathcal{X}(V)$ and $\bc\in E^7$ such that
$\wt(\bc)=4$ we have that
$$
(\bx\,|\,\bc)\in\mathcal{P}^\perp ~~\mbox{if and only if}~~ \bp(\bx,V) = \bc.
$$
\end{lemma}

\pr
Follows from the definition of the parity vector and orthogonality.
\qed

\medskip
We are ready to state now the main theorem.

\begin{theorem}\label{theo:1.1}
Let $v = 2^m-1 = 8\,u + 7 \geq 15$ and $k~=~u(u-1)/6$.

1).~ The number $M_{v,0}$ of all different Steiner triple systems
$S(v,3,2)$ of order $v$ and of rank $v-m$, orthogonal to the code $\mathcal{A}_m$, is equal to
$$
30^{u+1} \times 168^u\times 2^{3(u-m+3)}.
$$

2).~ The number $M_{v,1}$ of all different Steiner triple systems
$S(v,3,2)$ of order $v$ and of rank $v-m+1$, orthogonal to the code $\mathcal{A}_m$, is equal to
$$
210\times 6^u\times 840^u\times (8^u\times 16^{4k} - 7^u\times 2^{u-m+3})\times 2^{2(u-m+3)}.
$$

3).~ The number $M_{v,2}$ of all different Steiner triple systems
$S(v,3,2)$ of order $v$ and of rank $v-m+2$, orthogonal to the code $\mathcal{A}_m$, is equal to
$$
210\times [(4!\cdot 3!)^u \times 5040^u \times 2^{u-m+3}\times 576^{4k} - 48^u\times 2100^u\times 2^{2(u-m+3)}\times 16^{4k}].
$$
\end{theorem}

\pr
Suppose we are given a Steiner triple system $S_v = S(v,3,2)$ of order
$v=2^m-1=8u+7$ and of rank not greater than $2^m - m + 2$.
It implies that the dual code $S_v^{\perp}$ contains a $[v,m-3,d^{\perp}]$-code
$\mathcal{A}_m$ (see formula (2) of [4]).
Denote by $\mathcal{E}_u$ (so that $\mathcal{E}_u\subset E^{v}$) a linear space
generated by $u$ vectors of the form:
\begin{equation}\label{eq:3.2}
\begin{array}{ccccc}
(\bbo  &\bbz   &\ldots  &\bbz    &0000000)\\
(\bbz  &\bbo   &\ldots  &\bbz    &0000000)\\
\ldots &\ldots &\ldots  &\ldots  &\ldots\\
(\bbz  &\bbz   &\ldots  &\bbo    &0000000),
\end{array}
\end{equation}
where $\bbo = (1,1,1,1,1,1,1,1)$ and $\bbz = (0,0,0,0,0,0,0,0)$.

It is obvious that $\dim(\mathcal{E}_u)=u$ and that
$\mathcal{A}_m\subset \mathcal{E}_u$. Then there exist a linear subspace
$\mathcal{B}_{u-m+3}\subset \mathcal{E}_u$ of dimension $u-m+3$ such
that $\mathcal{E}_u=\mathcal{A}_m\oplus\mathcal{B}_{u-m+3}$.
Note that the space $\mathcal{E}_u$ is the kernel of the mapping $\psi$
of space $E^v$.

According to Construction $I(8)$ (see \cite{ZZ3}) the different Steiner triple systems
$S(v,3,2)$ which are orthogonal to the code $\mathcal{A}_m$,
are constructed either from the different $(3,2,64)_8$ codes $L$
or from $u$ different ordered partitions $V_1,V_2,\ldots,V_u$ of the set $W_2^8$
into parallel classes or from the different triple systems $S(7,3,2)$.
In particular, the different ordered partitions $V_1,V_2,\ldots,V_u$
correspond to the different sets $S^{(2,1)}$ (see $\S 4$ of \cite{ZZ3}).
Taking into account all possible orderings of partitions,
there are $(7!\cdot \Gamma_V)^u$ possible ways to construct the set $S^{(2,1)}$.
Note that the space $\mathcal{E}_u$ is by definition orthogonal to
any set $S^{(2,1)}$.

Suppose now, that we would like to construct a Steiner system whose
rank is less than $v-m+3$. It implies that in addition to $m-3$
orthogonal vectors given by the matrix (\ref{eq:3.1}), there exists one
(for the systems of rank $v-m+2$), two (for the systems of rank $v-m+1$)
or three (for the boolean systems of rank $v-m$) vectors extra vectors. These vectors
must be orthogonal to the system under construction and are located at a distance
$(v+1)/2$ from the words of the code $\mathcal{A}_m$. This condition is satisfied
if the extra orthogonal vector (or vectors) have weight $4$
in all $u+1$ blocks (the last block is of length $7$) (see [2, $\S 3$] and [4, $\S 4$]).

Since the non-zero coordinate positions of $S^{(3)}$ correspond to the last coordinate
positions of $\mathcal{E}_u$, the set $S^{(3)}$ is orthogonal to $\mathcal{E}_u$.
Note that the set $S^{(3)\perp}$ coincides with the parity matrix of the
default partition of $V^{(4)}$. We are going to consider the Steiner triple systems
whose set $S^{(3)}$ is given by (\ref{eq:S3}) and then multiply it by $30$.

Choose an arbitrary set $S^{(2,1)}$ generated by some partitions $V_1,V_2,\ldots,V_u$.
Suppose that there exists an orthogonal vector (to the set $S^{(2,1)}$)
$$
\bx=(\bx_1 ~|~ \ldots ~|~\bx_u ~|~ \bx_{u+1})\in E^v
$$
of weight $4$ in every of its $u+1$ blocks (i.e. $\wt(\bx_i)=4$, $i=1,\ldots,u+1$).
Recall Lemma \ref{lem:1.2}. Lemma \ref{lem:1.2} implies that for the
non-existence of an extra orthogonal vector $\bx$ it is enough that
among $u$ partitions $V_1,V_2,\ldots,V_u$ at least one of them is equivalent
to $V^{(5)}$ or to $V^{(6)}$. Thus
$$
\{V_1,V_2,\ldots,V_u\} ~\cap~ \bigl\{\Orb(V^{(5)})\,\cup\, \Orb(V^{(6)})\bigr\}\neq\emptyset.
$$

Note that the block $\bx_{u+1}\in E^7$ (the $(u+1)$-th block of $\bx$)
gets chosen according to the set $S^{(3)}$. Namely we can choose $\bx_{u+1}$
to be any of the seven non-zero vectors of a 3-dimensional space $S^{(3)\perp}$,
which is orthogonal to $S^{(3)}$. Since $S^{(2,1)}$ is orthogonal to $\bx$
every partition $V_i$ ($i=1,\ldots,u$) has respectively three $\bx_i$-even and four $\bx_i$-odd
parallel classes, i.e. we have
\begin{equation}\label{eq:3.2A}
\psi(\bx_i)\in \psi(\mathcal{X}(V_i))\subset E^8/{\la}\bbz,\bbo{\ra}.
\end{equation}
The sets $\mathcal{X}(V_i)$ are equivalent (up to the permutation of coordinates of $E^8$)
to one of the following sets: $\mathcal{X}(V^{(1)})$, $\mathcal{X}(V^{(2)})$, $\mathcal{X}(V^{(3)})$ or
$\mathcal{X}(V^{(4)})$.

For any index $i\in\{1,\ldots,u\}$, a set $P\subset E^8$ and vector $\be\in E^7$ define
$$
P\times_i\be = \{(\by_1 ~|~ \ldots ~|~\by_u ~|~ \be\in E^v ~|~ \by_i\in P,
\by_j=\bbz, j\in\{1,\ldots,u\}\backslash\{i\}\}.
$$
Using this notation, we have that
$$
S^{(2,1)}\, =\, \bigcup_{i=1}^u\,\bigcup_{j=1}^7 P_{i,j}\times_i\be_j,~\mbox{where}~ V_i= (P_{i,1},P_{i,1},\ldots,P_{i,7}).
$$

Fix any $i\in \{1,\ldots,u\}$. Consider all words of $S^{(2,1)}$ with nonzero
support at block $i$. According to Lemma \ref{lem:1.5}, the orthogonality condition to $\bx$ implies that
\begin{equation}\label{eq:Ortho1}
\bp_i(\bx_i,V_i) = \bx_{u+1}\in S^{(3)\perp},
\end{equation}
for $i=1,\ldots,u$. 
For any $\tau_i\in\mS_7$ consider $\tau_i\ast V_i$ (see (\ref{eq:perm})). According to Lemma \ref{lem:1.4}
$$
\bp_i = \bp_i(\bx_i,\tau\ast V_i) = \tau^{-1}(\bp_i(\bx_i,V_i)) =  \tau^{-1}(\bx_{u+1})\in S^{(3)\perp}.
$$
We conclude that $\tau_i^{-1} = \sigma_i\tau_i'$ where $\tau_i'\in\Stab(\bx_{u+1})$
and $\sigma_i$ is any fixed permutation such that
$$
\bp_i = \sigma_i(\bx_{u+1})\in S^{(3)\perp}.
$$
In other words $\tau_i$ can be any permutation from the coset $\sigma_i\Stab(\bx_{u+1})$.
In summary, given a fixed set of $u$ partitions $V_1,V_2,\ldots,V_u$, it is possible to construct $(4!\cdot 3!)^u$
different sets $S^{(2,1)}$ (using different orders of parallel classes $\tau_1\ast V_1,\ldots,\tau_u\ast V_u$)
which are orthogonal to the vector $\bx$.

One can see that if any vector $\bx$ satisfies conditions (\ref{eq:3.2A}),
then any vector $\bx' = \bx + \by$ for an arbitrary $\by\in\mathcal{E}_u$
satisfies conditions (\ref{eq:3.2A}) as well (since $\psi(\bx)=\psi(\bx')$).
And if $\by\in\mathcal{A}_m$, then ${\la}\mathcal{A}_m,\bx{\ra} = {\la}\mathcal{A}_m,\bx'{\ra}$,
i.e. in this case the vectors $\bx$ and $\bx'$ define the same orthogonal code,
and consequently the same Steiner triple system.
Thus in order to obtain different Steiner triple systems,
one should assume that $\by\in\mathcal{B}_{u-m+3}$.
So the number of ways to choose vector $\bx$, i.e. the blocks
$\bx_1,\ldots,\bx_u$ (where the block $\bx_{u+1}$ is fixed) is equal to
\begin{equation}\label{eq:3.2B}
|\mathcal{B}_{u-m+3}| \times |\psi(\mathcal{X}(V_1))| \times\ldots\times |\psi(\mathcal{X}(V_u))|.
\end{equation}

Suppose we have $u$ partitions $V_1,V_2,\ldots,V_u$ and assume
that we construct a set $S^{(2,1)}$ orthogonal to the vectors $\bx,\bx',\bx''\in E^v$
(which are linear independent in $E^v/\mathcal{E}_u$). Let
\begin{eqnarray*}\label{eq:bbb}
\bx   & = & (\bx_1 ~|~ \ldots ~|~\bx_u ~|~ \bx_{u+1}) \\
\bx'  & = & (\bx_1' ~|~ \ldots ~|~\bx_u' ~|~ \bx_{u+1}')\\
\bx'' & = & (\bx_1'' ~|~ \ldots ~|~\bx_u'' ~|~ \bx_{u+1}''),
\end{eqnarray*}
where $\wt(\bx_i)=\wt(\bx_i')=\wt(\bx_i'')=\wt(\bx_i+\bx_i')=\wt(\bx'_i+\bx_i'')=\wt(\bx_i+\bx_i'')=4$
for $i=1,\ldots,u+1$ and $\bx_{u+1}$, $\bx_{u+1}'$, $\bx_{u+1}''$ obviously satisfy
${\la}\bx_{u+1},\bx_{u+1}',\bx_{u+1}''{\ra}=S^{(3)\perp}$.

According to Lemma \ref{lem:1.2} this is possible when all partitions
$V_1,V_2,\ldots,V_u$ are equivalent to $V^{(4)}$ (i.e for any $i\in\{1,\ldots,u\}$,
$V_i\in\Orb(V^{(4)})$). For any $i$, consider three vectors $\bx_i,\bx_i',\bx_i''\in\mathcal{X}(V_i)$
linearly independent in $E^v/\mathcal{E}_u$ (i.e. such that
$\la\psi(\bx_i),\psi(\bx_i'),\psi(\bx_i'')\ra = \psi(\mathcal{X}(V_i))$).
Suppose that $\tau_i\in\mS_7$ and consider $\tau_i\ast V_i$ (see (\ref{eq:perm})). According to Lemma \ref{lem:1.4} and
due to the orthogonality condition (\ref{eq:Ortho1}), the parity vectors of $\tau_i\ast V_i$ are
equal
\begin{eqnarray*}
\la\bp_i,\bp_i',\bp_i''\ra
& = & \la\bp_i(\bx_i,\tau_i\ast V_i),\bp_i'(\bx_i',\tau_i\ast V_i),\bp_i''(\bx_i'',\tau_i\ast V_i)\ra\\
& = & \tau_i^{-1}(\la\bp_i(\bx_i,V_i),\bp_i'(\bx_i',V_i),\bp_i''(\bx_i'',V_i)\ra)\\
& = & \tau_i^{-1}(\la\bx_{u+1},\bx_{u+1}',\bx_{u+1}''\ra).
\end{eqnarray*}
We conclude that $\tau_i\in\Stab(S^{(3)\perp})$. Thus for any $i\in \{1,\ldots,u\}$ the parallel
classes of $V_i$ can be permuted by any permutation from $\Stab(S^{(3)\perp})$
(where $|\,\Stab(S^{(3)\perp})|=168$). So there exist $168^u$ different ways to construct
vectors $\bx,\bx',\bx''$ from the sets of blocks
$$
\{\bx_1,\bx_1',\bx_1''\},\ldots, \{\bx_u,\bx_u',\bx_u''\}
$$
and the corresponding permutations of parallel classes of
$\tau_1\ast V_1,\tau_2\ast V_2,\ldots,\tau_u\ast V_u$.

One can see that if any vectors $\bx,\bx',\bx''$ satisfy conditions (\ref{eq:3.2A}),
then the vectors $\bx + \by$, $\bx' + \by'$, $\bx''+ \by''$
for any $\by,\by',\by''\in\mathcal{E}_u$ satisfy conditions
(\ref{eq:3.2A}) as well. So the number of ways to choose vectors $\bx,\bx',\bx''$
and the corresponding permutations of the parallel classes of $V_1,V_2,\ldots,V_u$
(where the blocks $\bx_{u+1},\bx_{u+1}',\bx_{u+1}''$ are fixed) is equal to
\begin{equation}\label{eq:3.2BB}
168^u\times |\mathcal{B}_{u-m+3}|\times |\mathcal{B}_{u-m+3}| \times |\mathcal{B}_{u-m+3}| = 168^u\times 2^{3(u-m+3)}.
\end{equation}

Now suppose we have $u$ partitions $V_1,V_2,\ldots,V_u$ and assume
that we construct a set $S^{(2,1)}$ orthogonal to the vectors $\bx,\bx'\in E^v$
(which are linear independent in $E^v/\mathcal{E}_u$, i.e. $\psi(\bx)\neq\psi(\bx')$).
According to Lemma \ref{lem:1.2} this is possible when the set $S^{(2,1)}$ is constructed from partitions
$V_1,V_2,\ldots,V_u$ which are equivalent to $V^{(1)}$ or $V^{(4)}$.
Let
$$
\bx=(\bx_1 ~|~ \ldots ~|~\bx_u ~|~ \bx_{u+1}),~ \bx'=(\bx'_1 ~|~ \ldots ~|~\bx'_u ~|~ \bx'_{u+1}),
$$
where $\wt(\bx_i)=\wt(\bx_i')=\wt(\bx_i+\bx_i')=4$ for $i=1,\ldots,u+1$.
According to Lemma \ref{lem:dim2}, there are $7$ possibilities to choose
a pair of different vectors $\bx_{u+1},\bx_{u+1}'$ of $S^{(3)\perp}$.
Suppose that $\bx_i,\bx_i'\in\mathcal{X}(V_i)$.
Suppose that $\tau_i\in\mS_7$ and consider $\tau_i\ast V_i$. According to Lemma \ref{lem:1.4} and
due to the orthogonality condition (\ref{eq:Ortho1}), the parity vectors of $\tau_i\ast V_i$ are
equal to
\begin{eqnarray*}
\la\bp_i,\bp_i'\ra
& = & \la\bp_i(\bx_i,\tau\ast V_i),\bp_i'(\bx_i',\tau\ast V_i)\ra\\
& = & \tau^{-1}(\la\bp_i(\bx_i,V_i),\bp_i'(\bx_i',V_i)\ra)\\
& = & \tau^{-1}(\la\bx_{u+1},\bx_{u+1}'\ra).
\end{eqnarray*}
We conclude that $\tau_i^{-1} = \sigma_i\tau_i'$ where $\tau_i'\in\Stab(\la\bx_{u+1},\bx_{u+1}'\ra)$
and $\sigma_i$ is any fixed permutation such that
\begin{equation}\label{eq:3.2BC}
\la\bp_i,\bp_i'\ra = \sigma_i(\la\bx_{u+1},\bx_{u+1}'\ra).
\end{equation}
If $V_i$ is equivalent to $V^{(1)}$ the choice of $\bx_i,\bx_i'$ is unique.
According to Lemma \ref{lem:1.2A}, there exist $48$ possible permutations $\tau_i$ of the parallel classes of $V_i$
(all $\tau_i$'s belong to the coset of $\Stab(\la\bx_{u+1},\bx_{u+1}'\ra)$).

Suppose that $V_i$ is equivalent to $V^{(4)}$.
Then there exist seven two dimensional subspaces of a three
dimensional space ${\la}\psi(\mathcal{X}(V_i)){\ra}$.
Thus the number of ways to choose the pairs $\bx_i,\bx_i'$ and possible permutations $\tau_i$
of the parallel classes of $V_i$ is $7\times 48.$ Note that when
$$
\tau_i' \in \Stab(\bx_{u+1}) \,\cap\, \Stab(\bx_{u+1}'), ~\mbox{where}~ |\,\Stab(\bx_{u+1}) \,\cap\, \Stab(\bx_{u+1}')| = 8,
$$
all permutations $\tau_i\ast V_i$ correspond to the fixed pair of vectors $\bx,\bx'$.

Along with the vectors $\bx,\bx'$ one can also choose the vectors
$\bx+\by$ and $\bx'+\by'$ for any $\by,\by'\in\mathcal{B}_{u-m+3}$.
Taking into account Lemma \ref{lem:dim2} the number of different ways to choose $V_i$'s
and construct a pair of vectors $\bx,\bx'$ from the set of blocks
$$
\{\bx_1,\bx_1'\},\ldots, \{\bx_u,\bx_u'\}
$$
where the blocks $\bx_{u+1},\bx_{u+1}'$ are fixed, equals
\begin{equation}\label{eq:3.2D}
48^u\times (|\,\Orb(V^{(1)})| + 7\cdot |\,\Orb(V^{(4)})|)^u \times |\mathcal{B}_{u-m+3}| \times |\mathcal{B}_{u-m+3}|.
\end{equation}

Recall how the set $S^{(1,1,1)}$ (see [4, $\S 4$] for details),
which is orthogonal to one, two, or three extra vectors
(of weight $4$ in each block) gets constructed.
Let $S(u,3,2)$ be the boolean system uniquely defined
by the code $\mathcal{A}_{m}$. It is orthogonal
to the code obtained from $\mathcal{A}_{m}$ by erasing
the last seven coordinate positions (with zeroes)
and applying the map
$$
(0,0,0,0,0,0,0,0) \rightarrow 0, ~(1,1,1,1,1,1,1,1) \rightarrow 1.
$$
For a fixed vector $\bx=(\bx_1,\ldots,\bx_u)$ form $E^v$ (respectively for a fixed
pair or triple of orthogonal extra vectors linearly independent in $E^v/\mathcal{E}_u$)
and a fixed system
$S(u,3,2)$ of $k~=~u(u-1)/6$ vectors, every such vector
$\bc^{(j)}$ with support $\supp(\bc^{(j)}) = \{i_1,i_2,i_3\}$
generates the code $C(L_j;i_1,i_2,i_3)$ which is a subset
of a set $S^{(1,1,1)}$ of a new system $S(v,3,2)$.

According to Lemma \ref{lem:1.1}, the number of the different
codes $C(L_j;i_1,i_2,i_3)$ is equal to $576^4$ (respectively $16^4$ for a pair
of vectors and $1$ for a triple of extra orthogonal vectors linearly independent in $E^v/\mathcal{E}_u$).
For all $k$ vectors $\bc^{(j)}$ of the
system $S(u,3,2)$ there exist $576^{4k}$ (respectively $16^{4k}$ or $1$)
different codes $C(L_j;i_1,i_2,i_3)$. According to Lemma \ref{lem:1.1},
two different sets of codes $L_1,\ldots,L_k$ from the corresponding sets
$$
L_j \in \mathcal{L}((\bx_{i_1},\bx_{i_2},\bx_{i_3}))\subset\mathcal{L},~~~\{i_1,i_2,i_3\}=
\supp(\bc^{(j)}),~~~j=1,\ldots,k,
$$
define the different sets $S^{(1,1,1)}$ and consequently the different
resulting Steiner systems $S(v,3,2)$. Thus, to a fixed vector $\bx$
(respectively a pair or triple linearly independent in $E^v/\mathcal{E}_u$)
there corresponds $576^{4k}$ (respectively $16^{4k}$ or $1$) different sets
$S^{(1,1,1)}$ of rank not greater than $v-m+2$ (respectively $v-m+1$ or $v-m$).

All Steiner triple systems of rank $v-m$ have the orthogonal code
generated by $\mathcal{A}_{m}$ and three extra vectors from $E^v$
(linearly independent in $E^v/\mathcal{E}_u$):
\begin{eqnarray*}
\bx   & = & (\bx_1 ~|~ \ldots ~|~\bx_u ~|~ \bx_{u+1}) \\
\bx'  & = & (\bx_1' ~|~ \ldots ~|~\bx_u' ~|~ \bx_{u+1}')\\
\bx'' & = & (\bx_1'' ~|~ \ldots ~|~\bx_u'' ~|~ \bx_{u+1}''),
\end{eqnarray*}
where $\wt(\bx_i)=\wt(\bx_i')=\wt(\bx_i'')=\wt(\bx_i+\bx_i')=\wt(\bx'_i+\bx_i'')=\wt(\bx_i+\bx_i'')=4$
for $i=1,\ldots,u+1$.
According to Lemma \ref{lem:1.2} this is possible when the set $S^{(2,1)}$
is constructed from partitions $V_1,V_2,\ldots,V_u$ such that all of them
are equivalent to $V^{(4)}$ (i.e. $V_i\in\Orb(V^{(4)}))$.
The parallel classes of $V^{(4)}$ can be permuted by any permutation from $\Stab(S^{(3)\perp})$.
So there exist $30^u\times 168^u$ different ways to construct the set $S^{(2,1)}$.
There are $30$ possibilities to choose $S^{(3)}$ whereas three linearly independent blocks
$\bx_{u+1}$, $\bx_{u+1}'$ and $\bx_{u+1}''$ obviously satisfy
${\la}\bx_{u+1},\bx_{u+1}',\bx_{u+1}''{\ra}=S^{(3)\perp}$.
Since
$$
{\la}\psi(\bx_i),\psi(\bx_i'),\psi(\bx_i''){\ra} = {\la}\psi(\mathcal{X}(V_i)){\ra},
$$
there exist
$$
168^u\times |\mathcal{B}_{u-m+3}|^3 = 168^u\times 2^{3(u-m+3)}
$$
different ways to construct the triples of vectors $\bx,\bx',\bx''$ and the corresponding
partitions $\tau_i\ast V_1,\ldots,\tau_u\ast V_u$.
According to Lemma \ref{lem:1.1}, the set
$S^{(1,1,1)}$ in this case gets constructed uniquely. The first statement of this
Theorem follows.

Consider the Steiner triple systems of rank $v-m+1$ whose orthogonal code is generated by
$\mathcal{A}_{m}$ and a pair of vectors $\bx=(\bx_1 ~|~ \ldots ~|~\bx_u ~|~ \bx_{u+1})\in E^v$
and $\bx'=(\bx'_1 ~|~ \ldots ~|~\bx'_u ~|~ \bx'_{u+1})\in E^v$ (linearly independent in $E^v/\mathcal{E}_u$)
such that $\wt(\bx_i)=\wt(\bx_i')=\wt(\bx_i+\bx_i') = 4$.
As it was shown earlier, this is possible when the set $S^{(2,1)}$
is constructed from partitions $V_1,V_2,\ldots,V_u$ which are equivalent to $V^{(1)}$ or $V^{(4)}$.
The set $S^{(3)}$ can be arbitrary as well as the subspace
$\la\bx_{u+1},\bx'_{u+1}\ra\in S^{(3)\perp}$ where $\bx_{u+1}\neq\bx'_{u+1}$.

According to Lemma \ref{lem:1.1} there exist $16^{4k}$ different possibilities for the set $S^{(1,1,1)}$.
Using formula (\ref{eq:3.2D}) and taking into account that
$|\,\Orb(V^{(1)})| + 7\cdot |\,\Orb(V^{(4)})| = 840$,
the number of different Steiner triple systems is equal to
$$
30\cdot 7\times 48^u\times 2^{2(u-m+3)}\times 840^u\times 16^{4k}.
$$
The term $30\cdot 7$ is the number of different sets $S^{(3)}$ and the number of ways
to choose a subspace $\la\bx_{u+1},\bx_{u+1}'\ra\in S^{(3)\perp}$. Note that we have obtained the number
of all Steiner systems of rank not greater than $v-m+1$. However the systems of rank $v-m$
are counted multiple times. Recall that for the number $M_{v,0}$ of Steiner triple systems of rank $v-m$
(orthogonal to $\mathcal{A}_m$) we have
$$
M_{v,0} = 30^{u+1} \times 168^u\times 2^{3(u-m+3)}.
$$
As it was shown above, they are constructed from partitions $V_1,\ldots,V_u$ which are
equivalent to $V^{(4)}$. So, for every $i$, there are seven different $2$-dimensional subspaces
in $\psi(\mathcal{X}(V_i))$. Each of them can be extended to a $3$-dimensional space. So, every system
of rank $v-m$ is counted exactly $7\times 7^u$ times. Since $210\cdot 168 = 840\cdot 42$,
we obtain that the number $M_{v,1}$ of Steiner triple systems of rank $v-m+1$ orthogonal to $\mathcal{A}_m$
is equal to
$$
210\times 48^u\times 2^{2(u-m+3)}\times 840^u\times 16^{4k} - 210\times 210^u\times 168^u\times 2^{3(u-m+3)}
$$
$$
= 210\times 6^u\times 840^u\times 2^{2(u-m+3)}\times (8^u\times 16^{4k} - 7^u\times 2^{u-m+3}).
$$
The second statement of this Theorem follows.

Consider the Steiner triple systems of rank $v-m+2$ whose orthogonal code is generated by $\mathcal{A}_{m}$
and a vector $\bx=(\bx_1 ~|~ \ldots ~|~\bx_u ~|~ \bx_{u+1})\in E^v$ such that
$\wt(\bx_1)=\ldots=\wt(\bx_{u+1}) = 4$.
According to Lemma \ref{lem:1.2} this is possible when the set $S^{(2,1)}$ is constructed
from partitions $V_1,V_2,\ldots,V_u$ such that none of them are equivalent to $V^{(5)}$ or $V^{(6)}$.
For any $i$, there exist
$$
3\times|\,\Orb(V^{(1)})| + |\,\Orb(V^{(2)})|+ |\,\Orb(V^{(3)})| + 7\times|\,\Orb(V^{(4)})|
$$
$$
 = 3\times 630 + 420 +2520 + 7\times 30 = 5040,
$$
possibilities to choose vector $\bx_i$ and partition $V_i$ and
there exist $4!\cdot 3!$ permutations of its parallel classes.
According to Lemma \ref{lem:1.1} there are $576^{4k}$ different ways to construct the set $S^{(1,1,1)}$.
Thus the number of Steiner triple systems for this case is equal to
\begin{equation}\label{eq:3.3}
30\cdot 7\times (4!\cdot 3!)^u \times 5040^u \times 2^{u-m+3}\times 576^{4k}.
\end{equation}
The term $30\cdot 7$ is the number of different sets $S^{(3)}$ and blocks $\bx_{u+1}\in S^{(3)\perp}$.
Note that this expression is equal to the number
of different Steiner triple systems of rank at most $v-m+2$.

Now let's find how many these Steiner triple systems have rank $v-m$ or $v-m+1$.
The systems of rank $v-m$ occur when all partitions $V_1,\ldots,V_u$ are equivalent
to $V^{(4)}$. The spaces $\la\psi(\mathcal{X}(V_i))\ra$ have dimension three.
In this case the number of one dimensional
and two dimensional subspaces is seven. Thus, similar to the case of $v-m+1$,
every system of rank $v-m$ is counted exactly $7\times 7^u$ times.
However the multiplicities of the Steiner systems of rank $v-m+1$ occurring
in (\ref{eq:3.3}) are different for different systems. Let's compute them.

The systems of rank $v-m+1$ are constructed from partitions $V_1,\ldots,V_u$ which are equivalent
to either $V^{(1)}$ or $V^{(4)}$.
Suppose that $j$ partitions out of $V_1,\ldots,V_u$ are equivalent to $V^{(1)}$ whereas
the remaining $u-j$ are equivalent to $V^{(4)}$.

According to Lemma \ref{lem:1.1}, when the vectors $\bx,\bx'$ are chosen, the set $S^{(1,1,1)}$
can be constructed $16^{4k}$ different ways. When $j>0$, the resulting system is of rank $v-m+1$.
The number of such systems is equal to
$$
210\times 48^u\times {u\choose {j}}\times 630^j\times (7\cdot 30)^{u-j}\times 2^{2(u-m+3)}\times 16^{4k}.
$$
Now let's compute how many times any system of such type is counted (\ref{eq:3.3}).
If $V_i$ is equivalent to $V^{(1)}$, there are three non-zero vectors in $\la\psi(\mathcal{X}(V_i))\ra$,
whereas if $V_i$ is equivalent to $V^{(4)}$, the space $\la\psi(\mathcal{X}(V_i))\ra$ has seven non-zero
vectors and the same number of $2$-dimensional subspaces.
Thus, every such system is counted $3^j$ times in (\ref{eq:3.3}).

When $j=0$, the number of such systems of rank $v-m+1$ is equal to
$$
210\times 48^u\times (7\cdot 30)^{u}\times 2^{2(u-m+3)}\times 16^{4k} - 210\times 210^u\times 168^u\times 2^{3(u-m+3)}.
$$
Every such system is counted once, because the number of one and two dimensional subspaces is the same.
Since $3\cdot 630 + 7\cdot 30 = 2100$, we have
$$
\sum_{j=0}^{u} {u\choose j}\times (3\cdot 630)^{j} \times (7\cdot 30)^{u-j} = 2100^u.
$$
Thus the number of Steiner triple systems of rank $v-m+1$ with the corresponding multiplicities
which occur in (\ref{eq:3.3}) is equal to
$$
210\times 48^u\times 2100^u\times 2^{2(u-m+3)}\times 16^{4k} - 210\times 210^u\times 168^u\times 2^{3(u-m+3)}.
$$
Subtracting this number from (\ref{eq:3.3}) we obtain
\begin{eqnarray*}
M_{v,2}
& = & 210\times (4!\cdot 3!)^u \times 5040^u \times 2^{u-m+3}\times 576^{4k}\\
& - & 210\times[48^u\times 2100^u\times 2^{2(u-m+3)}\times 16^{4k} -  210^u\times 168^u\times 2^{3(u-m+3)}]\\
& - & 210\times 210^u\times 168^u\times 2^{3(u-m+3)}\\
& = & 210\times [(4!\cdot 3!)^u \times 5040^u \times 2^{u-m+3}\times 576^{4k}
- 48^u\times 2100^u\times 2^{2(u-m+3)}\times 16^{4k}].
\end{eqnarray*}
Thus the Theorem follows.
\qed


\nocite{}

\end{document}